\documentclass[a4paper,11pt]{amsart}

\usepackage{amsmath, amssymb, latexsym, amsthm, amsfonts, setspace, bm}
\usepackage[dvips]{color, graphicx}
\usepackage{float}

\theoremstyle{definition}
\newtheorem{definition}{Definition}[section]
\newtheorem{thm}[definition]{Theorem}

\newtheorem{lem}[definition]{Lemma}
\newtheorem{exmp}[definition]{Example}
\newtheorem{rem}[definition]{Remark}

\theoremstyle{definition}
\numberwithin{equation}{section}

\begin{document}
\title{A note on lower bounds for boxicity of graphs}
\author{Akira Kamibeppu}
\address{Department of Creative Engineering, National Institute of Technology, Kushiro College, Kushiro, Hokkaido 084-0916, Japan.}
\email{kamibeppu@kushiro-ct.ac.jp}

\date{}
\subjclass[2010]{05C62, 05C65, 05C72}
\keywords{hypergraph; boxicity; (co)interval graph; integer/linear program.}
\thanks{This work was supported by Grant-in-Aid for Young Scientists (B), No.25800091 (2013-2016).}

\begin{abstract}
The {\it boxicity} of a graph $G$ is the minimum non-negative integer $k$ such that $G$ can be isomorphic 
to the intersection graph of a family of boxes in Euclidean $k$-space, where a {\it box} in Euclidean $k$-space 
is the Cartesian product of $k$ closed intervals on the real line.\ 
In this note, we define the {\it fractional boxicity} of a graph as the optimum value of the linear relaxation 
of a covering problem with respect to boxicity, which gives a lower bound for its boxicity.\ 
We show that the fractional boxicity of a graph is at least the lower bounds for boxicity given by Adiga et al.\ in 2014.\  
We also present a natural lower bound for fractional boxicity of graphs.\ Moreover we discuss and focus on ``accuracy'' 
rather than ``simplicity'' of these lower bounds for boxicity as the next step in Adiga's work.\ 
\end{abstract}
\maketitle

\section{Introduction and Preliminaries}
A {\it box} in Euclidean $k$-space is the Cartesian product of $k$ closed intervals on the real line.\
The {\it boxicity} of a graph $G$, denoted by $\text{box}(G)$, is the minimum non-negative 
integer $k$ such that $G$ can be isomorphic to the intersection graph of a family of boxes in Euclidean $k$-space.\ 
The concept of boxicity of graphs was introduced by Roberts \cite{Ro69}.\ It has applications in some research fields, 
for example, a problem of niche overlap in ecology (see \cite{Ro76} for detail).\ So far many researchers have 
attempted to calculate or bound boxicity of graphs with specific structure, although we do not mention them here.\ 
Recently, Adiga et al.\ \cite{ACS14} presented a lower bound for the boxicity of a graph 
as in Lemma 1.1 below, which also gives some lower bounds under various conditions on graphs.\ 
Those lower bounds for boxicity in addition to the lower bound in Lemma 1.2 are relatively 
easy to estimate by examination, but there is an example of a graph whose boxicity 
cannot be determined by those lower bounds (see Example 2.4 and Remark 2.5).\ In what follows, 
the symbol $\overline{G}$ denotes the complement of a graph $G$.\ 
\begin{lem}[\cite{ACS14}, Lemma 3.1]
The inequality $\text{box}(G)\geq |E(\overline{G})|/|E(\overline{I_{\text{min}}})|$ holds for 
a non-complete graph $G$, where $I_{\text{min}}$ is an interval supergraph of $G$ with $V(I_{\text{min}})=V(G)$ 
and with the minimum number of edges among all such interval supergraphs of $G$.\ 
\end{lem}
\begin{lem}[\cite{CR83}, Lemma 3]
Let $G$ be a graph.\ Let $S_1=\{u_1, u_2, \ldots , u_n\}$ and $S_2=\{v_1,v_2, \ldots , v_n\}$ be disjoint 
subsets of $V(G)$ such that the only edges between $S_1$ and $S_2$ in $\overline{G}$ are the edges 
$u_iv_i$, where $i\in \{1,2,\ldots , n\}$.\ Then $\text{box}(G)\geq \lceil n/2 \rceil $ holds.\   
\end{lem}
The next step in Adiga's work is to discuss and focus on ``accuracy'' rather than ``simplicity'' of these 
lower bounds for boxicity.\ The purpose of this note is 
\begin{itemize}
\setlength{\itemsep}{0.1cm}
\item to review the lower bound in Lemma 1.1 for boxicity in the context of fractional graph theory, 
\item to introduce a fractional analogue of boxicity that will become a lower bound for boxicity, and  
\item to present a natural lower bound for our fractional analogue of boxicity, which works on calculation of 
boxicity of some graphs better than Lemma 1.1 and 1.2.\  
\end{itemize}
In this note, all graphs are finite, simple and undirected.\ We use $V(G)$ for the vertex set of 
a graph $G$ and $E(G)$ for the edge set of the graph $G$.\ These notations are also used for hypergraphs.\ 
A graph is said to be {\it interval} if its boxicity is at most 1.\ The cardinality of a set $X$ is denoted by $|X|$.\ 
A few concepts and results about (hyper)graphs are needed to present a fractional analogue of boxicity.\ 
A graph is said to be {\it cointerval} if its complement is an interval graph.\ A {\it cointerval edge covering} 
of a graph $G$ is a family $\mathcal{C}$ of cointerval subgraphs of $G$ such that each edge of $G$ is 
in some graph in $\mathcal{C}$.\ The following is a basic result on boxicity.\ 
\begin{thm}[\cite{CR83}, Theorem 3]Let $G$ be a graph.\ Then, $\text{box}(G)\leq k$ if and only if 
there exists a cointerval edge covering ${\mathcal C}$ of $\overline{G}$ with $|{\mathcal C}|=k$.\ 
Hence $\text{box}(G)=\min \{\,|\mathcal C| : \text{${\mathcal C}$ is a cointerval edge covering of $\overline{G}$\,}\}$.
\end{thm}

\section{Main Results}
Let ${\mathcal C}$ be a family of hyperedges of a hypergraph ${\mathcal H}$ and we write 
${\mathcal C}=\{X_1, \ldots , X_k\}$.\ The family ${\mathcal C}$ is a covering of ${\mathcal H}$ if 
$V({\mathcal H}) \subseteq X_1\cup \cdots \cup X_k$ holds.\ 
Our key idea for the definition of a fractional analogue of boxicity is in the way to define 
a hypergraph associated with a graph.\ For a graph $G$, we define the hypergraph ${\mathcal H}_G$ as follows:
\begin{align*}
V({\mathcal H}_G)&=E(\overline{G})\,\,\text{and}\,\,\\
E({\mathcal H}_G)&=\{E\subset E(\overline{G}) : \text{$E$ induces a cointerval subgraph of $\overline{G}$\,}\}. 
\end{align*}
Note that a covering of ${\mathcal H}_G$ corresponds to a cointerval edge covering of $\overline{G}$.\  
Hence the covering number of the hypergraph ${\mathcal H}_G$, the minimum cardinality of 
a covering of ${\mathcal H}_G$, is equal to the boxicity of $G$ by Theorem 1.3.\

For a graph $G$, let $e_i$ be an edge of $\overline{G}$ and $E_j$ a hyperedge of ${\mathcal H}_G$.\ 
Moreover, let $M_G$ be the incidence matrix of ${\mathcal H}_G$ whose rows are indexed by all edges of $\overline{G}$ 
and whose columns are indexed by all cointerval subgraphs of $\overline{G}$, that is, $i, j$-entry of $M_G$ 
is equal to 1 if $e_i\in E_j$, otherwise 0.\ Write $E({\mathcal H}_G)=\{E_1, \ldots , E_n\}$.\ 
Let ${\mathcal C}$ be a family of hyperedges in $E({\mathcal H}_G)$ and $\bm{x}=\,^t(x_1, x_2, \ldots , x_n) \in \{0,1\}^n$ 
the indicator vector of hyperedges in $E({\mathcal H}_G)$ that corresponds to the family ${\mathcal C}$, 
that is, $x_i$ is equal to 1 if $E_i\in {\mathcal C}$, otherwise 0.\ 
We see that ${\mathcal C}$ is a cointerval edge covering of $\overline{G}$ if and only if 
$M_G\bm{x}\geq 1$ (that is, each coordinate of $M_G\bm{x}$ is at least 1) holds.\ Hence the boxicity of a graph $G$ 
can be defined as the optimum value of the integer program
\begin{equation*} 
\begin{split}
\text{(IP)} \,\, & \text{minimize } ^t\bm{1x} \\
&\text{subject to}\,\, M_G\bm{x} \geq 1\,\, \text{and}\,\, \bm{x}\in \{0,1\}^n,  
\end{split}
\end{equation*}
that is,
\begin{equation*}
\text{box}(G)=\min \{{}^t\bm{1x} :  M_G\bm{x} \geq 1, \bm{x}\in \{0,1\}^n \}, 
\end{equation*}
where $\bm{1}$ is a vector of all ones.\ We relax the condition of the integer program (IP) and consider 
the linear program 
\begin{equation*} 
\begin{split}
\text{(LP)} \,\, & \text{minimize}\,\, ^t\bm{1x} \\
&\text{subject to}\,\, M_G\bm{x} \geq 1\,\, \text{and}\,\, \bm{x}\geq 0. 
\end{split}
\end{equation*}
We define the {\it fractional boxicity} of a graph $G$, denoted by $\text{box}_f(G)$, to be the optimum 
value of (LP), that is, 
\begin{equation*}
\text{box}_f(G)=\min \{{}^t\bm{1x} : M_G\bm{x} \geq 1, \bm{x}\geq 0 \}.
\end{equation*}
Hence $\text{box}_f(G)\leq \text{box}(G)$ holds for a graph $G$.\

By the way, in the theory of linear programming, we usually consider the dual program of (LP): 
\begin{equation*} 
\begin{split}
\text{(D)}\,\, &\text{maximize}\,\, ^t\bm{1y}\\ 
&\text{subject to}\,\, {}^tM_G\bm{y} \leq 1\,\, \text{and}\,\,\bm{y}\geq 0. 
\end{split}
\end{equation*}
It is well-known in the theory of linear programming that a linear program and its dual have the same optimum value.\ 
Hence we may consider the value of (D) instead of $\text{box}_f(G)$.\ 
We notice that a vector $\bm{y_*}$ of all $1/p$'s is a feasible solution of (D), where 
$p=\max_{E_i\in E({\mathcal H}_G)}|E_i|$.\ 
Hence, $\text{box}_f(G)\geq {}^t\bm{1y_*}=|E(\overline{G})|/p$.\
We note that this lower bound for fractional boxicity of graphs is identical to the lower bound for boxicity of graphs 
in Lemma 1.1.\

An automorphism of a hypergraph $\mathcal H$ is a bijection $\pi $ on $V(\mathcal H)$ such that $X\in E(\mathcal H)$ 
if and only if $\pi (X)\in E(\mathcal H)$.\
A hypergraph ${\mathcal H}$ is {\it vertex-transitive} ({\it edge-transitive}) 
if for every pair $(w_1, w_2)$ of vertices (hyperedges) there exists an automorphism $\pi$ of ${\mathcal H}$ 
such that $\pi (w_1)=w_2$ holds.\ The following theorem is derived from Proposition 1.3.4 in \cite{SU}.\ 

\begin{thm}
For a graph $G$, the inequalities
\begin{equation*}
\text{box}(G)\geq \text{box}_f(G)\geq \frac{|E(\overline{G})|}{\max_{E_i\in E({\mathcal H}_G)}|E_i| }
\end{equation*}
hold.\ Especially if $\overline{G}$ is edge-transitive, we have the equality
\begin{equation*}
\text{box}_f(G)=\frac{|E(\overline{G})|}{\max_{E_i\in E({\mathcal H}_G)}|E_i| }.
\end{equation*}
\begin{proof}
Note that the fractional boxicity of a graph $G$ is the same concept with the fractional covering number of 
the hypergraph ${\mathcal H}_G$.\ It is possible to show that the hypergraph ${\mathcal H}_G$ is 
vertex-transitive by the edge-transitivity of $\overline{G}$, so the above equality holds by 
Proposition 1.3.4 in \cite{SU}.\ The following Lemma 2.2 completes the proof of this theorem.\  
\end{proof}
\end{thm}
\begin{lem}
If $\overline{G}$ is edge-transitive for a graph $G$, the hypergraph ${\mathcal H}_G$ is vertex-transitive. 
\begin{proof}
For every pair of vertices $e_1, e_2\in V({\mathcal H}_G)=E(\overline{G})$, there exists an automorphism 
$\pi : V(\overline{G})\rightarrow V(\overline{G})$ such that $\pi (e_1)=e_2$ holds by our assumption.\  
We can check that $\pi $ induces a bijection $\overline{\pi }$ on $E(\overline{G})$ in a natural way: 
$\overline{\pi }(uv)=\pi (u)\pi (v)$.\ 
Moreover $E$ is in $E({\mathcal H}_G)$ if and only if $\overline{\pi }(E)$ is in $E({\mathcal H}_G)$ since $\pi $ 
and its inverse $\pi ^{-1}$ map a subgraph $H$ of $\overline{G}$ to the subgraph isomorphic to $H$.\ 
Hence $\overline{\pi }$ is the desired map.\ 
\end{proof}  
\end{lem}
The fractional boxicity of a graph $G$ is the same as the maximum value of $^t\bm{1y}$ under the conditions 
${}^tM_G\bm{y} \leq 1$ and $\bm{y}\geq 0$.\ We note that each entry of $\bm{y}$ is a weight of an edge 
of $\overline{G}$.\ The rows of ${}^tM_G$ are indexed by all cointerval subgraphs of $\overline{G}$, but we see that 
\begin{itemize}
\item an inequality 
in ${}^tM_G\bm{y} \leq 1$ corresponding to a non-maximal cointerval subgraph (on their edge sets) 
is superfluous since $\bm{y}\geq 0$.\ 
\end{itemize}
Hence we have only to focus on maximal cointerval subgraphs 
of $\overline{G}$ when we calculate $\text{box}_f(G)$.\ In what follows, $M_G$ always means the (reduced) 
incidence matrix of ${\mathcal H}_G$ whose columns are indexed by all maximal cointerval subgraphs of 
$\overline{G}$.\ 

We will reduce unnecessary restrictions further within the same conditions $M_G\bm{x}\geq 1$ and 
$\bm{x}\geq 0$.\ Let ${\mathcal E}$ ($\subset E({\mathcal H}_G)$) be the family of all maximal cointerval 
subgraphs of $\overline{G}$.\ Write ${\mathcal F}_e =\{E \in {\mathcal E} : e\in E \}$ for an edge 
$e\in E(\overline{G})$.\ An edge $e$ of $\overline{G}$ is said to be {\it fundamental} if ${\mathcal F}_e$ is 
minimal as subfamily of ${\mathcal E}$.\   
Let $E^*$ be the set of all fundamental edges of $\overline{G}$.\ We define two edges $e$ and $e'$ in $E^*$ 
to be equivalent, denoted by $e \sim e'$, if ${\mathcal F}_e = {\mathcal F}_{e'}$.\ We remark that 
\begin{itemize}
\setlength{\itemsep}{0.1cm}
\item an inequality in $M_G\bm{x}\geq 1$ corresponding to a non-fundamental edge of $\overline{G}$ is 
superfluous since $\bm{x}\geq 0$, and 
\item if $e \sim e'$ for $e$, $e' \in E^*$, 
the two inequalities in $M_G\bm{x}\geq 1$ which correspond to $e$ and $e'$ are the same inequalities.\ 
\end{itemize}
The inequality corresponding to an equivalence class $[e]$ means an inequality in $M_G\bm{x}\geq 1$ corresponding 
to a representative of $[e]$.\ It does not depend on the choice of representatives of $[e]$.\  
Let $M_G^*$ be the reduced incidence matrix of ${\mathcal H}_G$ whose rows are indexed by all equivalence 
classes in $E^*/{\sim }$ and whose columns are indexed by all maximal cointerval subgraphs of $\overline{G}$.\ 
We see that 
\begin{itemize}
\item $M_G\bm{x}\geq 1$ is equivalent to $M_G^*\bm{x}\geq 1$ under $\bm{x}\geq 0$.\  
\end{itemize}
Hence the fractional boxicity of a graph $G$ is the same as the optimum value of the linear program  
\begin{equation*} 
\begin{split}
\text{(LP)'} \,\, & \text{minimize}\,\, ^t\bm{1x} \\
&\text{subject to}\,\, M_G^*\bm{x} \geq 1\,\, \text{and}\,\, \bm{x}\geq 0. 
\end{split}
\end{equation*}

We consider a relaxation program of (LP)' and get a natural lower bound for fractional boxicity of graphs.\ 

\begin{thm} For a graph $G$, let $\{H_1, H_2, \ldots, H_l\}$ be the family of all maximal cointerval subgraphs of $\overline{G}$ and 
let $E^*/{\sim }=\{[e_1], [e_2], \ldots , [e_k]\}$.\ Let $a_i$ be the number of fundamental edges in $\{e_1, e_2, \ldots , e_k\}$ 
which are contained in $H_i$ for $i\in \{1,2,\ldots , l\}$.\ Then 
\begin{equation*}
\text{box}_f(G)\geq \frac{k}{a^*}
\end{equation*}
holds, where $a^*=\max\{a_1, a_2, \ldots ,a_l\}$.
\begin{proof}
Note that $\text{box}_f(G)=\min \{{}^t\bm{1x} : M_G^*\bm{x} \geq 1, \bm{x}\geq 0 \}$.\ 
Sum up all $k$ inequalities in $M_G^*\bm{x}\geq 1$, 
and then we obtain  
\begin{equation*} 
a^*({}^t\bm{1x})=a^*(x_1+x_2+\cdots +x_l)\geq a_1x_1+a_2x_2+\cdots +a_lx_l\geq k,    
\end{equation*}
where $\bm{x}=\,^t(x_1, x_2, \ldots , x_n)$.\ Hence ${}^t\bm{1x}\geq k/a^*$ holds, that is, $\text{box}_f(G)\geq k/a^*$.\ 
\end{proof}
\end{thm}


The fractional boxicity of a graph will measure its boxicity more accurately than the other lower bounds 
for boxicity given by Adiga et al.\ in 2014, although it is a difficult parameter to estimate by examination 
like the other fractional graph invariants. 

\begin{exmp}
We consider the graph $G_k$ whose complement is the graph in Figure 1 below (and is not edge-transitive), where $k\geq 4$.\ 
We will find all maximal cointerval subgraphs of $\overline{G_k}$ and prove $\text{box}_f(G_k)=k/2$.\ 
\begin{figure}[!ht]
\centering
\includegraphics[scale=0.72,clip]{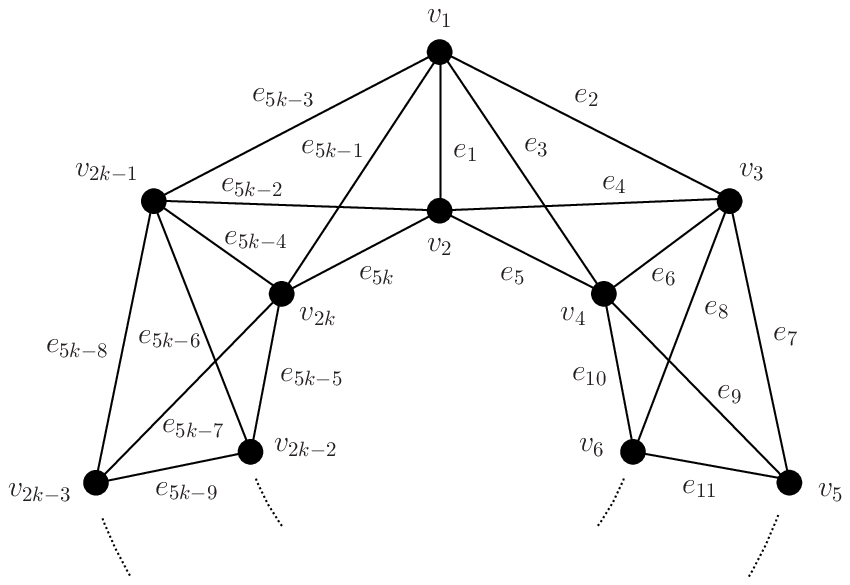}
\renewcommand{\baselinestretch}{1}
\caption[]{The complement of the graph $G_k$.}
\end{figure}\\
Let $H$ be a cointerval subgraph of $\overline{G_k}$.\ For example we see that 
\begin{itemize}
\setlength{\itemsep}{0.2cm}
\item[(1)]$H$ cannot have edges $e_{11}, e_{12}, \ldots , e_{5k-9}$ if $H$ has the edge $e_1$, and  
\item[(2)]$H$ cannot have edges $e_{16}, e_{17}, \ldots , e_{5k-9}$ if $H$ has at least one of $e_2, e_3, e_4$ or $e_5$.\ 
\end{itemize}
We will obtain similar statements to (1) or (2) if $H$ has an edge $e_i$, where $i\in \{6,7,\ldots , 5k\}$.\ \vspace{0.2cm}\\
\noindent {\bf Case I.} Assume that $H$ contains the edge $e_1$.\ If it has at least one of $e_7, e_8, e_9$ or $e_{10}$, 
we can find maximal cointerval graphs containing $H$ within the graph induced by $\{v_1, v_2, \ldots , v_6, v_{2k-1}, v_{2k}\}$, 
otherwise we can find them within the graph induced by $\{v_1, v_2, v_3, v_4, v_{2k-3}, v_{2k-2}, v_{2k-1}, v_{2k}\}$.\ \vspace{0.2cm}\\
\noindent {\bf Case II.} Assume that $H$ has at least one of $e_2, e_3, e_4$ or $e_5$.\ If it has at least one 
of $e_{12}, e_{13}, e_{14}$ or $e_{15}$, we can find maximal cointerval graphs containing $H$ within the graph 
induced by $\{v_1, v_2, \ldots , v_8\}$, otherwise we can find them within the graph 
induced by $\{v_1, v_2, \ldots , v_6, v_{2k-3}, v_{2k-2}, v_{2k-1}, v_{2k}\}$.\ \vspace{0.2cm}\\
As a result it is sufficient to find maximal cointerval subgraphs of the graph $H_*$ in Figure 2.\ 
Clearly, $H_*$ and $H_*-e$ (that is obtained from $H_*$ by deleting $e$) are not cointerval for any $e\in E(H_*)$.\ 
We will find three maximal cointerval subgraphs of $H_*$, but two of them can be extend to the graph isomorphic to 
the graph with heavy edges in Figure 2 on the original graph $\overline{G_k}$.\  
\begin{figure}[!h]
\centering
\includegraphics[scale=0.59,clip]{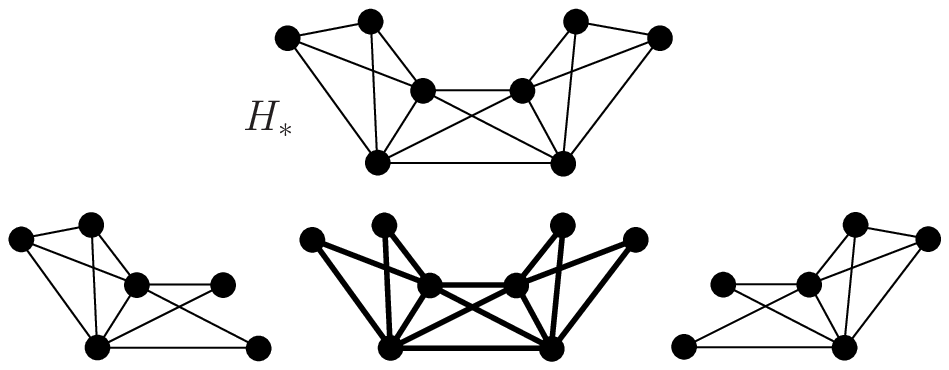}
\renewcommand{\baselinestretch}{1}
\caption[]{The graph $H_*$ and its maximal cointerval subgraphs.}
\end{figure}

We have $k$ maximal cointerval subgraphs of $\overline{G_k}$, each of which is isomorphic to the graph with 
heavy edges in Figure 2.\ Hence the optimum value of the following linear program becomes the fractional 
boxicity $\text{box}_f(G_k)$.\ 
\begin{align*}
\text{(D)}\,\,&\text{maximize} & & y_1+y_2+\cdots +y_{5k}\\
                     &\text{subject to} & & y_1+y_2 + \cdots +y_{10}+y_{5k-3}+y_{5k-2}+\cdots +y_{5k}\leq 1 \\
                     &                           & & y_{5i-13}+y_{5i-12}+\cdots +y_{5i}\leq 1 \hspace{0.5cm}(i\in \{3, 4, \ldots , k\})\\
                     &                           & & y_1+y_2+\cdots +y_5+y_{5k-8}+y_{5k-7}+\cdots +y_{5k}\leq 1\\
                     &                           & & y_j\geq 0\hspace{0.5cm}(j\in \{1,2,\ldots ,5k\})
\end{align*}
We consider the dual program of (D) and reduce superfluous inequalities so that we can obtain the following linear program:
\begin{align*}
\text{(LP)'}\,\,&\text{minimize}  & & x_1+x_2+\cdots +x_k\\
                     &\text{subject to} & & x_1+x_k\geq 1 \\
                     &                           & & x_i+x_{i+1}\geq 1 \hspace{0.5cm}(i\in \{1, 2, \ldots , k-1\})\\
                     &                           & & x_j\geq 0\hspace{0.5cm}(j\in \{1,2,\ldots ,k\}).
\end{align*}
Then $(x_1, x_2, \ldots, x_k)=(1/2, 1/2, \ldots , 1/2)$ is a feasible solution of (LP)', and hence $\text{box}_f(G_k)\leq k/2$.\ 

Let $E^*$ be the set of all fundamental edges of $\overline{G_k}$.\ It is easy to check that 
\begin{itemize}
\item $E^*=\{e_1, e_6, \ldots, e_{5k-4}\}$,  
\item for $e, e'\in E^*$, $e\ne e'$ implies $e\not\sim e'$, that is, \\
$E^*/{\sim } =\{[e_1], [e_6], \ldots, [e_{5k-4}]\}$, and 
\item $a^*=2$ since every maximal cointerval subgraph of $\overline{G_k}$ contains two fundamental edges 
in $\{e_1, e_6, \ldots, e_{5k-4}\}$.
\end{itemize}
By Theorem 2.3,  $\text{box}_f(G_k)\geq k/2$ holds, which implies our claim.\ \qed
\end{exmp}
\begin{rem}[$\text{box}_f$ vs.\ the lower bounds in Lemma 1.1 and 1.2]
It is easy to see that $\text{box}(G_k) \leq \lceil k/2\rceil $ holds for any $k$ by Theorem 1.3.\ 
Hence $\text{box}(G_k)=k/2$ holds since $\text{box}_f(G_k)=k/2$.\ The lower bounds 
in Lemma 1.1 and 1.2 do not work on the graph $G_k$ well, that is, they cannot determine the boxicity of $G_k$.\ 

We see that $\text{box}_f(G_k)>5k/14=|E(\overline{G_k})|/\max_{E_i\in E({\mathcal H}_{G_k})}|E_i|$ holds.\  
Let $m(G_k)$ be the maximum number of edges $a_ib_i$ of $G_k$ with the condition 
in Lemma 1.2 and let $M_k$ be a set of those edges.\  
For example if $e_1\in M_k$, any edge in $\{e_2, e_3, \ldots , e_{10}, e_{5k-8}, e_{5k-7}, \ldots , e_{5k}\}$ 
cannot be in $M_k$.\ If an edge $e\in \{e_2, e_3, e_4, e_5\}$ is in $M_k$, any edge in 
$\{e_1, e_2, \ldots , e_{11}, e_{5k-4}, e_{5k-3}, \ldots , e_{5k}\}\setminus \{e\}$ cannot be in $M_k$.\ 
It is not difficult to see that $m(G_k)\leq k/2$ holds.\ 
Then we have $\lceil m(G_k)/2 \rceil \leq \lceil k/4 \rceil <\text{box}_f(G_k)$.\
The difference between $\text{box}_f(G_k)$ and $|E(\overline{G_k})|/\max_{E_i\in E({\mathcal H}_{G_k})}|E_i|$ 
(or $\lceil m(G_k)/2 \rceil $) can be arbitrary large.\   
\end{rem}

\section{Further Observation}
Finally we remark another way to calculate the fractional boxicity of graphs.\ 
Let $s$ be a positive integer.\ The {\it $s$-fold boxicity} of a graph $G$, denoted by $\text{box}_s(G)$, 
is the minimum cardinality of a multiset $\{E_1, E_2, \ldots , E_k\}$ of cointerval subgraphs of 
$\overline{G}$ such that each edge of $\overline{G}$ is in at least $s$ cointerval subgraphs in the multiset.\ 
Note that $\text{box}_1(G)=\text{box}(G)$.\ 
Since the subadditivity $\text{box}_{s+t}(G)\leq \text{box}_s(G)+\text{box}_t(G)$ holds for a graph $G$ 
and $s,t\geq 1$,  the following limit exists and we have the following equality by Lemma 3.1:
\begin{equation*}
\lim _{s\rightarrow \infty} \frac{\text{box}_s(G)}{s}
=\inf \left\{\frac{\text{box}_s(G)}{s} : s\geq 1 \right \}.
\end{equation*}
\begin{lem}Let $\mathbb{Z^+}$ and $\mathbb{R^+}$ be the set of all nonnegative integers and 
the set of all nonnegative real numbers, respectively.\ If $g:\mathbb{Z^+}\rightarrow \mathbb{R^+}$ is 
subadditive, that is, $g(m+n)\leq g(m)+g(n)$ holds for any $m, n\in \mathbb{Z^+}$, the limit 
$\lim_{m\rightarrow \infty }g(m)/m$ exists and is equal to $\inf g(m)/m$. \qed
\end{lem} 
A basic result on the fractional covering numbers of hypergraphs guarantees 
$\text{box}_f(G)=\lim _{s\rightarrow \infty} \text{box}_s(G)/s$ (see Theorem 1.2.1 in \cite{SU}).\ 
Hence we may approach the study on the $s$-fold boxicity of graphs to calculate the fractional boxicity of graphs.\




\begin{thebibliography}{99}
\bibitem{ACS14}A.\,Adiga, L.\,S.\,Chandran and N.\,Sivadasan, Lower bounds for boxicity, Combinatorica 34 (2014) 631-655.
\bibitem{CR83}M.\,B.\,Cozzens and F.\,S.\,Roberts, Computing the boxicity of a graph by covering 
its complement by cointerval graphs, Discrete Appl.\,Math.\,6 (1983) 217-228.
\bibitem{Ro69}F.\,S.\,Roberts, On the boxicity and cubicity of a graph, in:\ Recent Progress in 
Combinatorics, Academic Press, New York (1969) 301-310.
\bibitem{Ro76}F.\,S.\,Roberts, Discrete mathematical models, with applications to social, biological,
and environmental problems, Prentice-Hall, New Jersey, 1976.
\bibitem{SU} E.\,R.\,Scheinerman and D.\,H.\,Ullman, Fractional Graph Theory, A rational approach to the theory of graphs, 
Dover Publications, Inc. Mineola, New York, 2011.  
\end{thebibliography}
\end{document}